\documentclass{amsart}

\usepackage[english]{babel}
\usepackage[utf8x]{inputenc}
\usepackage[T1]{fontenc}

\usepackage[a4paper,top=3cm,bottom=2cm,left=3cm,right=3cm,marginparwidth=1.75cm]{geometry}

\usepackage{amsmath}
\usepackage{amssymb,latexsym,amsmath}
\usepackage{graphicx}
\usepackage[colorinlistoftodos]{todonotes}
\usepackage[colorlinks=true, allcolors=blue]{hyperref}
\begin{document}

\title{An Alternate Method for Finding Particular Solutions of Nonhomogeneous Differential Equations}
\author{Ashot Djrbashian \and Milena Safaryan }

\date{}
\maketitle

\begin{abstract}
 In virtually every Ordinary Differential Equations textbook there are two main methods described for finding particular solutions of nonhomogeneous linear equations. They are the methods of Undetermined Coefficients and Variation of Parameters. Some books also present another method called Annihilator Method, but it is not very common by the reason of difficulty in application. In this short note we are presenting another method of finding particular solutions of nonhomogeneous linear equations. This approach  is different from the methods described above and, most importantly, it is universal in the sense that could be used in all circumstances. Moreover, in some very special cases, this approach can be used for equations with variable coefficients.\\



\end{abstract}

\noindent\textbf{1. INTRODUCTION.} Let us consider  an $n$th order linear nonhomogeneous differential equation with constant coefficients
\begin{equation}
a_ny^{(n)}+a_{n-1}y^{(n-1)}+\cdots+a_1y'+a_0y=q(t).
\end{equation}
 As it is well known, the general solution to this equation consists of the sum of the general solution of the homogeneous part and a particular solution of the equation (1). For this and all other facts used in this note we can refer, for example, to classical textbook [1]. While finding the general solution of the homogeneous part is reduced to solving an algebraic polynomial equation of order $n$ and is quite standard, finding the particular solution of  equation (1) sometimes creates  difficulties.  Two standard methods for this purpose are described in each textbook. These are the Method of Undetermined Coefficients and the Method of Variation of Parameters. Some textbooks also present another method, called the Annihilator Method, but it is relatively rare and not the most commonly used approach. In this paper we present another method that  reduces solution to repetitive solutions of first order equations. Moreover, we show that in some  special cases, this method can also be applied to find particular solutions of equations with variable coefficients. We will call this the  Differential Operator Approach (DOA) throughout this note.\\[10pt]

\noindent\textbf{2. DIFFERENTIAL OPERATOR APPROACH.} Let us begin with the simple case of the second-order equation: 

\begin{equation}
y'' + a y' + b y = q(t),
\end{equation}
Using the standard notation $Df=\frac{df}{dt}=f'$ we can rewrite this equation in the form 

\begin{equation}
\left( D^2 + aD + bI \right)y = q(t),
\end{equation}
where $I$ denotes the identity operator. The characteristic polynomial of the given differential equation is:

\[
p(r)=r^2 + ar + b,
\]
which can be factored as follows: $(r - r_1)(r - r_2)$,
where $ r_1 $ and  $r_2 $ are the roots of the polynomial. Consequently, equation (3) can be written as
$$(D-r_1I)(D-r_2I)y=q,$$
where the product of differential operators is understood as the composition of operators. 

Expanding the last equation, we will get  
\begin{equation}
    (D-r_1I)(D-r_2I)y=D((D-r_2I)y)-r_1((D-r_2I)y)=q.
\end{equation}    

Let  us denote $(D-r_2I)y=\phi.$ With this notation  equation (4) becomes 
$$(D-r_1I)\phi=\phi'-r_1\phi=q.$$
This means that finding a particular solution to the nonhomogeneous equation is reduced  to solving two first-order constant coefficient linear equations consequently.

Therefore, we begin by solving the "outer" equation first: $\phi'-r_1\phi=q(t).$ By the standard method, the solution is given as (we omit integration constants in the process of solution because we are  looking for particular solutions)
\begin{equation}
    \phi(t)=e^{r_1t}\int e^{-r_1t}q(t)dt.
\end{equation}
Using this solution, we can solve the "inner" equation $y'-r_2y=\phi(t).$ It is given, as above, using the formula 
\begin{equation}
 y=e^{r_2t}\int e^{-r_2t}\phi(t)dt.  
\end{equation}
Combining solutions (5) and (6) we obtain the particular solution of equations (2)--(3)

\begin{equation}
y_p = e^{r_2 t}  \int e^{(r_1-r_2) t} \left( \int e^{-r_1 t} q(t) \, dt \right) dt.
\end{equation}
From the demonstration of the solution of the second-order equation it is clear that the same approach  works for equations of any order. Moreover, the solution is independent from the type of "outside force" function $q$. Hence, this approach can substitute both the method of undetermined coefficients and the method of variation of parameters.
Below we demonstrate this approach using several examples.\\[6pt]

\noindent\textbf{3.EXAMPLES.} 1. \textbf{The case of distinct real roots.} To solve the differential equation 
\begin{equation}
y'' + 5y' + 6y = e^t\cos t
\end{equation}
we first solve the corresponding characteristic equation $r^2+5r+6=0$ and find $r_1=-2, r_2=-3$. 
Using the formula above, we see that a particular solution is given by 
$$ y_p(t)= e^{-2t}\int e^{2t}e^{-3t}\left(\int e^{3t}(e^t\cos t)dt\right)dt.$$
Calculating the inside integral gives (here and everywhere else below we disregard the integration constants)
$$\int e^{3t}(e^t\cos t)dt=\int e^{4t}\cos t dt=\frac{e^t}{17}(\sin t+4\cos t).$$

Continuing the process, we now have to calculate the "outer" integral
$$y_p(t)=e^{-2t}\int e^{2t}\frac{e^t}{17}(\sin t+4\cos t)dt=\frac{e^t}{170}(11\cos t+7\sin t).$$

Substituting into the original equation (8) shows that this is indeed a solution to the nonhomogeneous equation. \\[6pt] 


\noindent 2. \textbf{The case of repeated roots.} We now consider the case of repeated roots. The equation
\begin{equation}
y''-4y'+4y=t^3e^{2t}
\end{equation}
has a repeated root $r=2$, hence $e^{2t}$ and $te^{2t}$ are the solutions of the corresponding homogeneous equation. If we try to find a particular solution to (9) using the Method of Undetermined Coefficients, the expected form of the function would be $t^2(At^3+Bt^2+Ct+D)e^{2t}$. Therefore, to find a particular solution using this approach, we potentially need to solve a system of four linear algebraic equations with four unknowns A, B, C, and D. Using the DOA  instead, we have:
$$
\phi'-2\phi=t^3e^{2t}
$$
and the solution of this equation would be
$$
\phi(t)=e^{2t}\int e^{-2t}t^3e^{2t}dt=e^{2t}t^4/4.
$$
Next, we have to solve the equation $y'-2y=\phi$ using the same exact formula:
$$
y(t)=e^{2t}\int e^{-2t}e^{2t}t^4/4 dt=\frac{e^{2t}}{4}\int t^4 dt=\frac{1}{20}t^5e^{2t}
$$
and it's done! \\[6pt]
3. \textbf{The case of complex roots.} Let us now consider the case of complex roots of the characteristic equation. The equation
\begin{equation}
y''-2y'+5y=\sin t
\end{equation}
has the roots  $1\pm2i$ and the differential operator factorization will be
$$(D-(-1+2i)I)(D-(-1-2i)I)y=\sin t.$$
The standard method of solving this equation involves the process of excluding complex values and considering only the real solutions, but here, for the sake of simplifying integration, we will use  complex exponential solutions. 
$$y(t)=c_1e^{(-1+2i)t}+c_2e^{(-1-2i)t}.$$
Denoting, as before, $(D-(-1-2i)I)y=\phi$, we will have the first equation:
$$\phi'-(1+2i)\phi=\frac{e^{it}-e^{-it}}{2i},$$
where we used Euler's formula for the sine function. Solving this equation yields
$$\phi(t)=\frac{e^{(1+2i)t}}{2i}\int\left[e^{-(1+i)t}-e^{-(1+3i)t}\right]dt$$
$$=\frac{1}{2i}\left(\frac{e^{it}}{-1-i}-\frac{e^{-it}}{-1-3i}\right).$$
The second equation 
$$y'-(1-2i)y=\frac{1}{2i}\left(\frac{e^{it}}{-1-i}-\frac{e^{-it}}{-1-3i}\right)$$
has the solution
$$y=\frac{e^{(1-2i)t}}{2i}\int \left[\frac{e^{(-1+3i)t}}{-1-i}-\frac{e^{(-1+i)t}}{-1-3i}\right]dt$$
$$=\frac{e^{(1-2i)t}}{2i}\left[\frac{e^{(-1+3i)t}}{(-1-i)((-1+3i)}-\frac{e^{(-1+i)t}}{(-1-3i)(-1+i)}\right]$$
$$=\frac{1}{2i}\left[\frac{e^{it}}{4-2i}-\frac{e^{-it}}{4+2i}\right].$$
The denominators inside the brackets are complex conjugates. Therefore, by  combining the two fractions and further simplifying the expressions, we  arrive at the following final answer:
$$y_p(t)=\frac{1}{10}\cos t+\frac{1}{5}\sin t.$$
\\[6pt]
4. \textbf{The case when the undetermined coefficients method does not apply.}
Consider the example 
\begin{equation}
    y''+4y'+4y=e^{-2t}\ln t.
\end{equation}
For this kind of right-hand side functions only the variation of parameters method can be applied. Also, pay attention, that the characteristic equation has a repeated root $r=-2.$ Now, by re-writing in operator form, we obtain
\[
(D^2+4D+4I)y=D(D+2I)y+2(D+2I)y=e^{-2t}\ln t.
\]

Denoting $\phi=(D+2I)y$ we arrive to the first-order equation
\[
\phi '+2\phi=e^{-2t}\ln t
\]
which has the solution $\phi=e^{-2t}(t\ln t-t).$ Now, on the last step we solve the equation 
\[
y'+2y=e^{-2t}(t\ln t-t)
\]
which gives us the final solution
\[
y_p=\frac{1}{4}t^2e^{-2t}(2\ln t-3).
\]

It is also important  to notice, that in this particular case (and any other case when the variation of parameters method is used) we completely avoided the necessity of calculating the Wronskian of the homogeneous equation. Additionally, there is no need to memorize the two formulas that lead to the final solution in cases where that poarticular method is used. \\[6pt]
5. \textbf{The case of variable coefficients.} From the demonstration at the beginning of Section 2, it is clear that it is impossible to factor equations into a product of two (or more) lower-degree equations  when coefficients are variable. Moreover, the notion of the characteristic equation and characteristic roots do not apply here. As is well known, and described in any standard Differential Equations textbook, the only more or less productive method of solving variable coefficients equations is the Method of Power Series. Obviously, the Cauchy-Euler equation is the one lucky exception. Because of the difficulty in finding these kinds of solutions, the case of nonhomogeneous equations is rarely  mentioned in textbooks. Below we  show that in some limited special cases the DOA approach described above works here as well.\\
Let us consider the expression 
\begin{equation}
(D+ax^nI)(D+bx^mI)
\end{equation}
and see for which real values of $a, b, n, m$ this expression might take the form $D^2-cx^kI.$ Expanding the expression in (12) (using  the product rule) we obtain
\[
D^2+(ax^n+bx^m)D+(anx^{n-1}+bmx^{m-1}+abx^{n+m})I
\]
and the necessary conditions are  $n=m, b=-a.$ This means that the equations of the form $y''-a^2x^{2n}y=q$ can be factored into as 
\[
(D+ax^nI)(D-ax^nI)y=q(x).
\]
From here we can proceed exactly as in the case of constant variable coefficients and have $D(D-ax^nI)y+ax^n(D-ax^nI)y=q.$ Notation $\phi=(D-a^nxI)y$ reduces our equation to a first-order linear equation 
\[
\phi '+ax^n\phi=q(x).
\]
A particular solution to this equation is given by the formula
\[
\phi=\exp\left\{-\frac{ax^{n+1}}{n+1}\right\}\int \exp\left\{\frac{ax^{n+1}}{n+1}\right\}q(x)dx.
\]
In the second, and last step of the solution, we need to take this function and solve the second first order equation $y'-ax^ny=\phi.$ Once again, the solution is standard and is given by the following formula 
\[
y_p=\exp\left\{\frac{ax^{n+1}}{n+1}\right\}\int \exp\left\{-\frac{ax^{n+1}}{n+1}\right\}\phi(x)dx.
\]
\\[6pt]

\noindent\textbf{Concluding remarks}.
As demonstrated, the approach using linear differential operators effectively eliminates the need for guessing solutions and provides a systematic method to solve differential equations analytically. That said, it is quite possible that this method is not completely new. However, all our efforts to find a reliable reference where this was done before, were unsuccessful. We will appreciate if anybody could provide that kind of information.\\[6pt]


\noindent\textbf{REFERENCES}

\noindent\hrulefill

1. E. A. Coddington, \emph{An Introduction to Ordinary Differential Equations},  Dover Publications, New York, NY, 1989.
\\[6pt]

\it{Dr. Ashot Djrbashian, Division of Mathematics, Glendale Community College, 1500 N. Verdugo Rd., Glendale, CA 91208, ashotd@glendale.edu}

Milena Safaryan, Division of Mathematics, Glendale Community College, 1500 N. Verdugo Rd.,Glendale, CA 91208  msafary632@student.glendale.edu

\end{document}